\documentclass[final]{raa}          
\usepackage{graphicx,times}             

\begin{document}

   \title{RECURSIVE HARMONIC ANALYSIS FOR COMPUTATIONAL
OF HANSEN COEFFICIENTS
}

   \volnopage{Vol.0 (200x) No.0, 000--000}      
   \setcounter{page}{1}           

   \author{Mohamed Adel Sharaf
      \inst{1}
   \and Hadia Hassan Selim
      \inst{2}\mailto{hassanselim@hotmail.com}
      }

   \institute{ Astronomy Department, Faculty of Science, King Abdul Aziz University, Jeddah, Saudi Arabia\\
             \email{Sharaf$ \_ $adel@hotmail.com}
        \and
            Astronomy Department, National Research Institute of Astronomy and Geophysics,Helwan, Egypt\\
             \email{hassanselim@hotmail.com}
          }

   \date{Received~~2001 month day; accepted~~2001~~month day}

   \abstract{
This paper reports on a simple pure numerical method developed for
computing Hansen coefficients by using recursive harmonic analysis
technique. The precision criteria of the computations are very
satisfactory and provide materials for computing Hansen's and
Hansen's like expansions, also to check the accuracy of some
existing algorithms.
   }

   \authorrunning{M. A.Sharaf \& H. H. Selim }            
   \titlerunning{RECURSIVE HARMONIC ANALYSIS }  

   \maketitle

%
%
\section{Introduction}           
\label{sect:intro}
 Hansen coefficient( Cefola \cite{ce77})  is an
important class of functions which frequently occur in many
branches of Celestial Mechanics
        such as planetary theory(Newcomb \cite{ne95}) and artificial satellite motion (Allan\cite{al67}; Hughes\cite{hu77}).
        Moreover , there are extensive forms of  Hansen
        like expansions (Klioner et. al.~\cite{kl98} ; Sharaf \cite{sh85}, \cite{sh86}) which play important roles in the
         expansion theories  of elliptic motion.\\
        Giacalia (\cite{gi76}) noted that Hansen's coefficients appears in satellite theory in expression of the disturbing
        function due to the primary and due to the presence of a third body and they are usually called Eccentricity Functions.
        He derived recurrence relation for these functions and their derivatives, as they appear in the evaluation of geopotential
        and third body perturbations of an artificial satellite. Also in \cite{gia87}, he proved Hansen's coefficients for Fourier
        series in terms of the mean anomaly correspond to a rotation of the orbital plane proportional to the eccentricity of the
        orbit. They are given in terms of Bessel functions and generalized associated Legendre functions which arise
        through the transformation of spherical harmonics under rotation. In \cite{hug81}, Hughes computed tables of analytical expressions
        for the Hansen coefficients $\textsc{x}_{o}^{n,\pm m}(e)$
        and $\textsc{x}_{o}^{-(n+1),\pm m}(e)$ when $1\leq n\leq 30$
        and $0\leq m\leq n$. In \cite{bra90}, Branham derived a recursive
        calculation of Hansen coefficients which are used in expansions of elliptic motion by three methods: Tisserand's
        method, Von Zeipel-Andoyer method with explicit representation of the polynomials required to compute the Hansen
        coefficients and von Zeipel-Andoyer method with the value of the polynomials calculated recursively. Vakhidov
        (\cite{vak00})
        studied in detail efficient approximations of Hansen coefficients using polynomials in terms of the eccentricity.
        He and Zhang (\cite{mia90}) used Hansen coefficients to compute general perturbations of the asteroids of
        Flora group due to Jupiter. Breiter et.al (\cite{bre04}) show that most of the theory of Hansen coefficients remains valid for
        $\textsc{x}_{k}^{\gamma j}$, when $\gamma$ is a real
        number, also, the generalized coefficients can be applied in a variety of perturbed problems that involve some
        drag effects. Sadov (\cite{sad08}) deals analytically with the properties of Hansen's coefficients in the theory of
        elliptic motion considered as functions of the parameter
        $\eta=\sqrt{1-e^{2}}$ where $e$ is the eccentricity.\\
        In the present paper, we develop a simple pure numerical method for computing Hansen coefficients by using recursive
        harmonic analysis technique. The precision criteria of the computations are very satisfactory. The importance of
        the method is that it not only provides materials for computing Hansen's and also Hansen's like expansions but also,
        it can be used due to its simplicity and accuracy, to check the accuracy of the different algorithms already existing.


\section{Basic Formulations}
\label{sect:Bas}
\subsection{Properties of Least -Squares}

Let $y$ be represented by the general linear expression of the
form $\sum_{i=1}^{L}c_{i}\phi(x)$ where $\phi's$ are linear
independent functions of $x$. Let $ \textbf{c}$ be the vector of
the exact values of the $ c's$ coefficients and
$\widehat{\textbf{c}}$
 be the least -squares estimators of
 $\textbf{c}$ obtained from the solution of the normal equations $
 \textbf{G}\widehat{\textbf{c}}=\textbf{b}$ . The coefficient matrix
$ \textbf{G}(L\times L)$ is symmetric positive definite, that is,
all its eigenvalues $ V_{i}; i= 1,2,..,L$ are positive. Let
$\textsc{E C}(z)$ denote the expectation of $z$ and $\sigma^{2}$
the variance of the fit, defined as

\begin{equation}
\sigma^{2}=q_{n}/(N-L)
\end{equation}

where

\begin{equation}
q_{n}=(\textbf{y}- \Phi^{T}\widehat{\textbf{c}})^{T}(\textbf{y}-
\Phi^{T}\widehat{\textbf{c}})
\end{equation}

$N$ is the number of observations, $\textbf{y}$ is a vector with
elements $y_{k}$ and $\Phi(L\times N)$  has elements  $
\varphi_{ik}=\varphi_{i}(x_{k})$ . The transpose of a vector or a
matrix is indicated by the superscript $ 'T'$.\\
 According to the least- squares criterion, it could be shown that(Sharaf
et.al.\cite{sh00})
\begin{itemize}
 \item The estimators $\widehat{\textbf{c}}$ given by the least- squares method give the minimum of $q_{n}$.
 \item The estimators $\widehat{\textbf{c}}$ of the coefficients $\textbf{c}$, obtained by the method of least-squares, are unbiased; $i.e. EC(\widehat{\textbf{c}})=\textbf{c}$
 \item The variance-covariance matrix $Var(\widehat{\textbf{c}})$ of the unbiased estimators $\widehat{\textbf{c}}$ is given by

 \begin{equation}
 Var(\widehat{\textbf{c}})=\sigma^{2}\textbf{G}^{-1},
 \end{equation}

 where $\textbf{G}^{-1}$ is the inverse of $\textbf{G}$.
 \item The average squared distance between $\textbf{c}$ and $\widehat{\textbf{c}}$ is

\begin{equation}
EC(D^{2})=\sigma^{2}\sum_{1}^{L}\frac{1}{V_{i}}.
\end{equation}

\end{itemize}
\subsection{Harmonic Analysis of a Periodic Function}
Let it be required to find a sum

\begin{equation}
a_{o}+\sum_{j=1}^{s}a_{j}\cos jx+\sum_{j=1}^{s}b_{j}\sin jx
\end{equation}

 which furnishes the best possible representation of a function
 $u(x)$, when we are given that $u(x)$ takes the values
 $u_{o},u_{1},...,u_{i-1}$ when $x$ takes $x_{o},x_{1},...,
 x_{i-1}$ respectively , $m$ being some number greater than $2s$. The problem is to determine the
 $(2s+1)$ constants,  $a_{o},a_{j}$ and $b_{j} ; j=1,2,...,s$ so as to make the
 expression (5) takes, as nearly as possible, the  $l$ values
 $u_{o},u_{1},...,u_{i-1}$ when $x$ takes the values
 $x_{o},x_{1},...,x_{i-1}$. To do so we shall make use of the method of least squares and we get

 \begin{eqnarray}
 \frac{1}{2}a_{o}\eta_{oi}+\sum_{j=1}^{s}a_{j}\eta_{ij}+\sum_{j=1}^{s}b_{j}\beta_{ij}=d_{i},
 i=0,1,...,s ; \nonumber\\
 \frac{1}{2}a_{o}\beta_{oq}+\sum_{j=1}^{s}a_{j}\beta_{qj}+\sum_{j=1}^{s}b_{j}\gamma_{qj}=c_{q},
 q=1,2,...,s;
 \end{eqnarray}

where

\begin{eqnarray}
\eta_{ij}&=&\eta_{ji}=\sum_{k=0}^{i-1}\cos ix_{k}, i=0,1,...,s
,j=0,1,...s ;\nonumber\\
\beta_{qj}&=&\sum_{k=0}^{i-1}\cos jx_{k}
\sin qx_{k} , j=0,1,...,s , q=1,2,...,s ;\nonumber\\
\gamma_{qj}&=&\gamma_{jq}=\sum_{k=0}^{i-1}\sin qx_{k}\sin jx_{k} ,
q=1,2,...,s  ,  j=1,2,...,s;\nonumber\\
d_{i}&=&\sum_{k=0}^{i-1}u_{k}\cos ix_{k} ,i=0,1,...,s;\nonumber\\
c_{q}&=&\sum_{k=0}^{i-1}u_{k}\sin qx_{k}  , q=1,2,...,s.
\end{eqnarray}

Equations (7) are the normal equations. These equations represent
a set of linear equations in $(2s+1)$ unknowns $a's$ and $b's$
coefficients and could be solved by any of the methods adopted for
linear systems. However, the coefficient matrix of this set could
be reduced to a diagonal one by certain choice of the arguments
$x_{k}$ and in this case the $a's$ and $b's$ are determined
exactly and the problem is known as harmonic analysis.\\
In the method of harmonic analysis, the arguments $x_{k}$ take the
special values;

\begin{equation}
0,\frac{2\pi}{l}, 2.\frac{2\pi}{l},
3.\frac{2\pi}{l},...,(l-1).\frac{2\pi}{l}.
\end{equation}

For these values, the $\eta's , \beta's$ and $\gamma's$ of
Equations (7) become:\\
For $i=j\neq 0: \eta_{ij}=\gamma_{ji}=\frac{1}{2}l ; \beta_{ij}=0
.$\\
For $i\neq j : \eta_{ij}=\gamma_{ij}=\beta_{ij}=0$\\
Consequently the $a's$ and $b's$ coefficients could then be
computed exactly from

\begin{eqnarray}
a_{j}=\frac{\mu}{l}\sum_{k=0}^{i-1}u_{k}\cos j.\frac{2\pi}{l}k
,j=0,1,..,s ; \nonumber\\
b_{q}=\frac{2}{l}\sum_{k=0}^{i-1}u_{k}\sin q.\frac{2\pi}{l}k
,q=1,2,...,s.
\end{eqnarray}

where $\mu=1$ if $j=0$ ; $\mu=2$ if $j>0$ .

\subsection{Hansen Coefficients}

Consider elliptic motion expansions of $(r/a)^{n} \cos mv$ and
$(r/a)^{n}\sin mv$ in terms of the mean anomaly $M$  that is,

\begin{eqnarray}
(\frac{r}{a})^{n}\cos mv=\sum_{k=0}A_{k}^{n,m}\cos k M;\nonumber\\
(\frac{r}{a})^{n}\sin mv=\sum_{k=1}B_{k}^{n,m}\sin k M .
\end{eqnarray}

where $a$ is the semi major -axis , $r$ the radial distance, $n$
is a positive or negative integer ,while $m$ is positive integer
and $v$ the true anomaly in elliptic motion . The A's and B's
coefficients called Hansen's coefficients ,are functions of the
eccentricity $e$.\\
The relations between the eccentric anomaly  $E$ and the anomalies
$M,v$ are given for elliptic motion as follows:
\begin{itemize}
\item The relation between E and M is well know Kepler's equation

\begin{eqnarray}
M=E- e \sin E .
\end{eqnarray}

\item The fundamental relations between $v$ and $E$ in an elliptic orbit are

\begin{equation}
\tan  \frac{v}{2}=\sqrt{\frac{1+e}{1-e}}\tan \frac{E}{2}
\end{equation}

These equations are the most useful relations between $v(E)$ and
$E(v)$ , since $\frac{v}{2}$ and $\frac{E}{2}$ are always in the
same quadrant. There is a possibility of numerical trouble when
Equation (12) is used with angles that are near $\pm
\frac{\pi}{2}$ as the two tangents become infinite. In order to
avoid this difficulty, Broucke and Cefola \cite{br73} established
the formula

\begin{equation}
\tan \frac{1}{2}(v-E)=\frac{\beta \sin E}{1-\beta \cos E},
\end{equation}

where

\begin{equation}
\beta= \frac{1-\sqrt{1-e^{2}}}{e}=\frac{e}{1+\sqrt{1-e^{2}}}
\end{equation}

Equation (13) is free of numerical trouble, no matter what the
values of the angles are. Moreover, it can be easily used because
the angle $(v-E)/2$ is always less than $\frac{\pi}{2}$ for all
elliptic orbits.\\

\item Finally the relation between $r$ and $E$ is

\begin{equation}
(\frac{r}{a})= 1-e \cos E
\end{equation}

\end{itemize}

\section{Computational Developments}
\subsection{Practical Computations of the $a's$ and $b's$ Coefficients}
The $a's$ and $b's$ coefficients of Equations (9) could be
computed efficiently (Ralston \& Rabinowitz\cite{ra78})from

\begin{equation}
a_{j}=\frac{\mu}{l}\{u_{o}+F_{1,j} \cos \frac{2\pi}{l}j-F_{2,j}\};
j=0,1,...,s ,
\end{equation}

\begin{equation}
b_{q}=\frac{2}{l}F_{1,q} \sin \frac{2\pi}{l}q ; q=1,2,...,s.
\end{equation}

where, for any $j$ the $F's$ are computed recursively from

\begin{equation}
F_{k,j}= u_{k}+2 \cos x_{j} F_{k+1,j}-F_{k+2,j}
\end{equation}

by using the initial conditions $F_{i,j}=F_{i+1,j}=0$, starting
with $k=l-1$ to compute successively $F_{i-1,j}, F_{i-2,j}, ...,
F_{1,j}.$

\subsection{Error Estimates}
\begin{itemize}
\item The variance of the fit (Equation (1)) is given by

\begin{equation}
\sigma^{2}=\frac{\delta^{2}}{l-L}
\end{equation}

where the sum of the squares of the residuals $\delta^{2}$  is
given as (Ralston \& Rabinowitz\cite{ra78})

\begin{equation}
\delta^{2}=\sum_{i=0}^{i-1}u_{i}^{2}-\frac{1}{2}[2a_{o}^{2}+\sum_{j=1}^{s}(a_{j}^{2}+b_{j}^{2})]
\end{equation}

Clearly both $\sigma$ and $\delta$  depend on the number $s$ of
the $a's$ and $b's$ coefficients. If the precision is measured by
probable error $PE$, then

\begin{equation}
PE=0.6745 \sigma
\end{equation}

\item Since the coefficient matrix $\textbf{G}$ of the harmonic analysis is diagonal with
elements of the same value $l/2$ , then according to Equation (3)
the standard error of each of the $a's$ and $b's$ coefficients is

\begin{equation}
\sigma_{coeff}=\sigma\sqrt{\frac{2}{l}}
\end{equation}

The corresponding probable error for each coefficient is

\begin{equation}
PE_{coeff}=0.6745\sigma_{coeff}
\end{equation}

\item The average squared distance between the exact and the least -squares values
 (Equation (9)) is given according to Equation (4) as

 \begin{equation}
 Q=CE (D^{2})=\frac{2s}{l}\sigma^{2}
 \end{equation}

 \end{itemize}

\subsection{Choosing the Number of the Coefficients}

In practice ,since we do not know $s$ , we would evaluate $a's$
and $b's$ coefficients for $s=1,2,…$ , then compute $\delta^{2}$
(Equation(19)), and continue as long as $\delta^{2}$  decreases
significantly(within a given tolerance $Tol$) with increasing $s$
.
\subsection{The Special Values}
The special values of the left hand sides of Equation (10) are
computed as follows:

\begin{enumerate}
  \item $M_{i}=\frac{2\pi i}{l}$; $i=0,1,...,l-1$ .
  \item For each $M$ solve Kepler's equation (Equation(11)) by
Newton-Raphson iterative method (or any other method). Let $E_{o}$
be an initial approximation of $E$; define for $k=0,1,2,...$

   \begin{eqnarray*}
   E_{k+1}=E_{k}-\frac{E_{k}-e \sin E_{k}-M}{1-e \sin E_{k}}.
   \end{eqnarray*}

   Each $E_{k+1}$ should approximate $E$  more closely than $E_{k}$  . For the initial approximation
   $E_{o}$ use the value (Battin\cite{ba99})

   \begin{eqnarray*}
   E_{o}=M +\frac{e \sin M}{1- \sin (M+e) + \sin M}.
   \end{eqnarray*}

The above procedure is terminated if the following conditions are
satisfied \\
$\varepsilon_{2}\leq \varepsilon_{1}$ and $\mid
H(E_{i+1)}\mid\leq100 \varepsilon_{1}, \varepsilon_{2}=\mid
\frac{E_{i+1}-E_{i}}{E_{i+1}}\mid$ if $\mid E_{i+1}\mid > 1 ;
\varepsilon_{2}=\mid E_{i+1}-E_{i}\mid$ if $\mid E_{i+1}\mid<1 ,$
where $\varepsilon_{1}$ is a given tolerance and $H(E)=M-E-e \sin
E$.
  \item For each $E$ compute $v$ using Equation (14) and $(\frac{r}{a})^{n}$   from Equation (10)
  \item For each $v$ compute $\cos (mv)$.
  \item Finally, find the product of the values of $(\frac{r}{a})^{n}$(of step 3)and $\cos(mv)$ (of step 4).
\end{enumerate}

\subsection{Numerical Results}

The above computational developments are applied for calculating
Hansen's coefficients of Equation (10) with input constants as
$l=100 , Tol=10^{-6}$ and $\varepsilon_{1}=10^{-8}$. The numerical
results are listed in Tables $I$ to Table $VI$ for different
values of $n$, $m$ and different eccentricities of some members of
the solar system. In these tables $\delta_{A}^{2}(\delta_{B}^{2})$
represents the sum of the squares of the residuals of Equation
(19) for $A's ( B's )$ coefficients, $\sigma_{coeff.
A}(\sigma_{coeff. B})$ represents the common standard  error of
Equation( 21) for $A's ( B's )$ coefficients, finally
$Q_{A}(Q_{B})$ represents the average squared distance between the
exact and least- squares values of Equation (23) for $A's ( B's )$
coefficients.

\section{Conclusion}
In concluding the present paper, pure numerical method is
developed for computing Hansen coefficients by using recursive
harmonic analysis technique. The precision criteria which are: the
variance of the fit, the standard errors of the coefficients and
the average squared distance between the exact and least squares
values, are all very satisfactory. The  method is not only provide
materials for computing Hansen's and also Hansen's  like
expansions but also can be used to check the accuracy of the
different algorithms that already exist.

\begin{table}[h!]
\begin{center}
\bf{\sc Table I:}
\sc{Hansen Coefficients for The Planet Earth :\\
 $e=.016708617, n=-3, m=6$} \\
\begin{tabular}{cll}\\
\hline\noalign{\smallskip}
   $k$       &      $A_{k}$   &    $B_{k}$ \\\hline
0  &  $-2.80505\times10^{-16}$  &$\mbox{}$ \\
1  &  $-1.69927\times10^{-10}$  &    $-1.69927\times10^{-10}$    \\
2 & $1.31491\times10^{-7}$ & $1.31491\times10^{-7}$\\
3 &  $-0.0000259261$  &   $-0.0000259261$     \\
4 &  $0.00209013$ &   $0.00209013$     \\
5 & $-0.0749101$       & $-0.0749101$     \\
6 & $0.99039$   &  $0.99039$   \\
7 & $0.124591$   &  $0.124591$   \\
8 & $0.00917108$    &  $0.00917108$   \\
9 & $0.000516607$   &$0.000516607$    \\
10 & $0.0000246565$   &  $0.0000246565$  \\
11 & $1.05004\times10^{-6}$ & $1.05004\times10^{-6}$\\
\hline
$d_{A}^{2}=8.52652\times10^{-14}$ && $d_{B}^{2}=6.39488\times10^{-14}$\\
$s_{\infty ee.A}=4.37729\times10^{-9}$ && $s_{\infty ee.B}=3.79085\times10^{-9}$\\
$Q_{A}=2.10768\times10^{-16}$&& $Q_{B}=1.8076\times10^{-16}$\\
\hline
\end{tabular}
\end{center}
\end{table}

\begin{table}[h!]
\begin{center}
\bf{\sc Table II:}
\sc{Hansen Coefficients for The Planet Pluto : \\
$e=0.249050, n=-3, m=6$} \\
\begin{tabular}{cll}\\
\hline\noalign{\smallskip}
   $k$       &      $A_{k}$   &    $B_{k}$ \\\hline
0  &  $0.0508079$  & $\mbox{}$ \\
1  &  $-0.325005$  &  $-0.319177$    \\
2 & $0.969155$ & $0.969203$\\
3 &  $-1.34716$  &   $-1.34716$     \\
4 &  $0.536896$ &   $0.536896$     \\
5 & $0.248699$       & $0.248699$     \\
6 & $0.0765758$   &  $0.0765758$   \\
7 & $0.0211903$   &  $0.0211903$   \\
8 & $0.0056458$    &  $0.0056458$   \\
9 & $0.00148234$   &$0.00148234$    \\
10 & $0.000386931$   &  $0.000386931$  \\
11 & $0.000100722$ & $0.000100722$\\
12& $0.0000261571$ & $0.0000261571$ \\
13& $6.76933\times10^{-6}$ & $6.76933\times10^{-6}$ \\\hline
$d_{A}^{2}=1.62174\times10^{-10}$ && $d_{B}^{2}=1.6226\times10^{-10}$\\
$s_{\infty ee.A}=1.93084\times10^{-7}$ & & $s_{\infty ee.B}=1.93135\times10^{-7}$\\
$Q_{A}=4.84659\times10^{-13}$ && $Q_{B}=4.84914\times10^{-13}$\\
\hline
\end{tabular}
\end{center}
\end{table}

\begin{table}[h!]
\begin{center}
\bf{\sc Table III:}
\sc{Hansen Coefficients for The Asteroid Ceres:\\ $e=0.078, n=8, m=2$} \\
\begin{tabular}{cll}\\
\hline\noalign{\smallskip}
   $k$       &      $A_{k}$   &    $B_{k}$ \\\hline
0  &  $0.0854431$  & $\mbox{}$ \\
1  &  $-0.492936$  &  $-0.479094$    \\
2 & $1.08609$ & $1.08564$\\
3 &  $-0.157994$  &   $-0.157993$     \\
4 &  $0.00140598$ &   $0.00140603$     \\
5 & $0.000192711$       & $0.000192714$     \\
6 & $0.0000113508$   &  $0.0000113509$   \\
7 & $6.01265\times10^{-7}$   &  $6.01269\times10^{-7}$   \\
\hline
$d_{A}^{2}=4.26326\times10^{-14}$ && $d_{B}^{2}=4.26326\times10^{-14}$\\
$s_{\infty ee.A}=3.02792\times10^{-9}$& & $s_{\infty ee.B}=3.02792\times10^{-9}$\\
$Q_{A}=6.41781\times10^{-17}$ && $Q_{B}=6.41781\times10^{-17}$\\
\hline
\end{tabular}
\end{center}
\end{table}

\begin{table}[h!]
\begin{center}
\bf{\sc Table IV:}
\sc{Hansen Coefficients for The Asteroid Sekhmet :\\
$e=0.296, n=-1, m=5$} \\
\begin{tabular}{cll}\\
\hline\noalign{\smallskip}
   $k$       &      $A_{k}$   &    $B_{k}$ \\\hline
0  &  $-0.0000795273$  & $\mbox{}$ \\
1  &  $0.00983893$  &  $0.00983416$    \\
2 & $-0.114213$ & $-0.114214$\\
3 &  $0.431088$  &   $0.431088$     \\
4 &  $-0.482649$ &   $-0.482649$     \\
5 & $-0.260258$       & $-0.260258$     \\
6 & $0.191795$   &  $0.191795$   \\
7 & $0.410314$   &  $0.410314$   \\
8 & $0.411965$    &  $0.411965$   \\
9 & $0.318733$   &$0.318733$    \\
10 & $0.213837$   &  $0.213837$  \\
11 & $0.130854$ & $0.130854$\\
12 & $0.0750457$ & $0.0750457$ \\
13 & $0.0410069$ & $0.0410069$ \\
14 & $0.021582 $ & $0.021582$ \\
15 & $0.0110231$ & $0.0110231$ \\
16 & $0.00549385$ & $0.00549385$ \\
17 & $0.00268279$ & $0.00268279$\\
18 & $0.00128767$ & $0.00128767$ \\
19 & $ 0.000608982$ & $0.000608982$\\
20 & $0.000284349$ & $0.000284349$\\
21 & $0.000131294$ & $0.000131294$\\
22 & $0.0000600289$ & $0.0000600289$\\
23 & $0.0000272069$ & $0.0000272069$\\
24 & $0.0000122351$ & $0.0000122351$\\
25 & $5.46368\times10^{-6}$ & $5.46368\times10^{-6}$\\\hline
$d_{A}^{2}=3.64729\times10^{-10}$& & $d_{B}^{2}=3.64665\times10^{-10}$\\
$s_{\infty ee.A}=3.11867\times10^{-7}$& & $s_{\infty ee.B}=3.1184\times10^{-7}$\\
$Q_{A}=2.43152\times10^{-12}$ && $Q_{B}=2.4311\times10^{-12}$\\
\hline
\end{tabular}
\end{center}
\end{table}

\begin{table}[h!]
\begin{center}
\bf{\sc Table V:}
\sc{Hansen Coefficients for The Comet Wild2 : \\
$e=0.541, n=3, m=2$} \\
\begin{tabular}{ccc}\\
\hline\noalign{\smallskip}
   $k$       &      $A_{k}$   &    $B_{k}$ \\\hline
0  &  $-0.187235$  & $\mbox{}$ \\
1  &  $0.954443$  &  $0.943797$    \\
2 & $-1.75935$ & $-1.75982$\\
3 &  $1.04451$  & $1.04448$     \\
4 &  $0.271626$ & $0.271625$     \\
5 & $-0.101094$  & $-0.101092$     \\
6 & $-0.158792$   &  $-0.158791$   \\
7 & $-0.112669$   &  $-0.112669$   \\
8 & $-0.05443$    &  $-0.0544298$   \\
9 & $-0.0112087$   &$-0.0112086$    \\
10 & $0.0142725$   &  $0.0142726$  \\
11 & $0.026263$ & $0.026263$\\
12 & $0.0296931$ & $0.0296931$ \\
13 & $0.0283862$ & $0.0283862$ \\
14 & $0.0248699 $ & $0.0248699$ \\
15 & $0.0206484$ & $0.0206484$ \\
16 & $0.0165297$ & $0.0165297$ \\
17 & $0.0128893$ & $0.0128893$\\
18 & $0.00985376$ & $0.00985376$ \\
19 & $ 0.00741822$ & $ 0.00741822$\\
20 & $0.005551671$ & $0.005551671$\\
21 & $0.00406198$ & $0.00406198$\\
22 & $0.00296636$ & $0.00296636$\\
23 & $0.00215139$ & $0.00215139$\\
24 & $0.00155123$ & $0.00155123$\\
25 & $0.00111291$ & $0.00111291$\\
26 & $0.000795003$ & $0.000795003$\\
27 & $0.000565771$ & $0.000565771$\\
28 & $0.000401309$ & $0.000401309$\\
29 & $0.000283824$ &  $0.000283824$\\
30 & $0.000200214$ & $0.000200214$\\
31 & $0.000140908$ & $0.000140908$\\
32 & $0.0000989644$ & $0.0000989644$\\
33 & $0.0000693761$ & $0.0000693761$\\
34 & $0.000048552$ &  $0.000048552$\\
35 & $0.0000339265$ & $0.0000339265$\\
36 & $0.0000236736$ & $0.0000236736$\\
37 & $0.0000164982$ & $0.0000164982$\\
38 & $0.0000114842$ & $0.0000114842$\\
39 & $7.9854\times10^{-6}$ & $7.9854\times10^{-6}$\\\hline
$d_{A}^{2}=1.18559\times10^{-8}$ & &$d_{B}^{2}=1.18562\times10^{-8}$\\
$s_{\infty ee.A}=4.05228\times10^{-7}$ && $s_{\infty ee.B}=4.05232\times10^{-7}$\\
$Q_{A}=6.40418\times10^{-12}$ && $Q_{B}=6.4043\times10^{-12}$\\
\hline
\end{tabular}
\end{center}
\end{table}

\begin{table}[h!]
\begin{center}
\bf{\sc Table VI:}
\sc{Hansen Coefficients for The Comet Lexell :\\
 $e=0.786, n=8, m=4$} \\
\begin{tabular}{cll}\\
\hline\noalign{\smallskip}
   $k$       &      $A_{k}$   &    $B_{k}$ \\\hline
0  &  $28.4068$  & $\mbox{}$ \\
1  &  $-47.0631$  &  $-25.693$    \\
2 & $23.9162$ & $21.1464$\\
3 &  $-4.70405$  & $-4.84203$     \\
4 &  $-0.605464$ & $-0.619241$     \\
5 & $-0.0262285$  & $-0.0283248$     \\
6 & $0.0293643$   &  $0.0289445$   \\
7 & $0.0217703$   &  $0.0216686$   \\
8 & $0.0121887$    &  $0.0121605$   \\
9 & $0.00637396$   &$0.00636531$    \\
10 & $0.00325957$   &  $0.00325672$  \\
11 & $0.00164622$ & $0.00164524$\\
12 & $0.000817397$ & $0.000817052$ \\
13 & $0.000392683$ & $0.000392564$ \\
14 & $0.000176143 $ & $0.000176105$ \\
15 & $0.0000672271$ & $0.0000672181$ \\
16 & $0.0000140575$ & $0.0000140581$ \\
17 & $-0.0000103329$ & $-0.0000103297$\\
18 & $-0.0000200483$ & $-0.000020045$ \\
19 & $ -0.0000224845$ & $-0.0000224818$\\
20 & $-0.0000215081$ & $-0.0000215061$\\
21 & $-0.0000191168$ & $-0.0000191154$\\
22 & $-0.0000163167$ & $-0.0000163158$\\
23 & $-0.0000135893$ & $-0.0000135887$\\
24 & $-0.0000111418$ & $-0.0000111414$\\
25 & $-9.0409\times10^{-6}$ & $-9.04071\times10^{-6}$\\
\hline
$d_{A}^{2}=7.42148\times10^{-9}$ & &$d_{B}^{2}=7.33417\times10^{-9}$\\
$s_{\infty ee.A}=1.40679\times10^{-6}$ && $s_{\infty ee.B}=1.39849\times10^{-6}$\\
$Q_{A}=4.94765\times10^{-11}$ && $Q_{B}=4.88944\times10^{-11}$\\
\hline
\end{tabular}
\end{center}
\end{table}

\label{lastpage}

\end{document}